\documentclass{article}
\usepackage[utf8]{inputenc}
\usepackage{amssymb}
\usepackage{amsthm}
\usepackage{tikz}
\usepackage[hyphens]{url}
\usepackage[hyphenbreaks]{breakurl}
\usepackage{hyperref}
\usepackage[style=numeric,sorting=nty]{biblatex}

\newtheorem{theorem}{Theorem}
\newtheorem{corollary}{Corollary}
\newtheorem{lemma}{Lemma}

\newtheorem{conj}{Conjecture}
\newtheorem{propn}{Proposition}

\title{On a problem posed by Bjorn Poonen}
\author{James Rawson}

\begin{document}

\maketitle

\begin{abstract}
Bjorn Poonen asked whether there exists a polynomial giving a surjection $\mathbb{Z} \times \mathbb{Z} \to \mathbb{N}$. We answer this question in the negative, conditional on a conjecture of Vojta. More precisely, we show that if such a function exists, there is a family of open surfaces with dense integral points despite the surfaces being of log general type.
\end{abstract}
In 2009, Bjorn Poonen asked on MathOverflow \cite{poonen_q}, if there is a polynomial, $F \in \mathbb{Q}[x, y]$, such that $F(\mathbb{Z}, \mathbb{Z}) = \mathbb{N}$. This question is similar to two other questions, are there polynomial injections $\mathbb{Q} \times \mathbb{Q} \to \mathbb{Q}$, and are there polynomial bijections between the same? The answers to these are, respectively: yes, conditional on the Bombieri-Lang conjecture \cite{poonen_inj}, and no, also assuming the Bombieri-Lang conjecture \cite{bresciani}. In a similar vein, we answer this question in the negative, assuming an integral generalisation of the Bombieri-Lang conjecture.

After increasing the number of variables, it is classically known that such surjections exist: for $\mathbb{Z}^3 \to \mathbb{N}$, this can be done via the sum of 3 triangular numbers; and for $\mathbb{Z}^4 \to \mathbb{N}$ it can also be done by the sum of 4 squares. 

Throughout, we fix a polynomial $F \in \mathbb{Q}[x, y]$, and define the surface $X_n \subset \mathbb{A}^3$ by $z^n = F(x, y)$. We start by showing the following constraint on $F$ (see Section 1). 
\begin{theorem}
Suppose integral points are not dense on $X_n$ for some odd integer $n$, then $F(\mathbb{Z}, \mathbb{Z}) \neq \mathbb{N}$.
\label{fewpts}
\end{theorem}

Non-density of integral points is known to hold in certain situations, for example if $X_n$ is an abelian variety \cite{silverman}. In our case, we rely on the following conjecture of Vojta \cite{vojta}.
\begin{conj}[Vojta]
Let $V$ be a smooth variety, $\tilde{V}$ a smooth projective closure of $V$, $D$ the divisor supported on $\tilde{V} \setminus V$ and $K_{\tilde{V}}$ be the canonical divisor of $\tilde{V}$. If $\overline{K}_V = K_{\tilde{V}} + D$ is a big divisor (i.e. $V$ is of log general type), then integral points are not dense on $V$.
\end{conj}

This is an integral generalisation of the following conjecture of Bombieri and Lang.

\begin{conj}[Bombieri-Lang]
Let $V$ be a smooth projective variety and $K_{V}$ be the canonical divisor of $V$. If $K_V$ is a big divisor ($V$ is of general type), then rational points are not dense on $V$.
\end{conj}

A consequence of Theorem 1 and Vojta's conjecture is the following (proven in Section 2).
\begin{propn}
Assume Vojta's conjecture, then if $F$ satisfies $F(\mathbb{Z}, \mathbb{Z}) = \mathbb{N}$, then $X_n$ is never of log general type. Moreover, $\mathbb{A}^2 \setminus \{F(x, y) = 0\}$ is an open variety of logarithmic Kodaira dimension $-\infty$, 0 or 1.
\end{propn}

Polynomials in $\overline{\mathbb{Q}}[x, y]$ satisfying this last condition have been classified by Aoki \cite{aoki}. They show that $F(x, y) = a$ is a rational curve for all but finitely many $a$ ($F$ is generically rational), and either the normalisation is missing at most 2 points, or there is an auxilliary polynomial satisfying certain strong hypotheses. In Section 3, we show that, in each case, such polynomials in $\mathbb{Q}[x, y]$ either do not represent all of $\mathbb{N}$ or represent arbitrarily large negative integers.

\section{Integral Points Are Not Dense on $X_n$}
In this section, we will prove Theorem~\ref{fewpts}. Throughout this section, $n$ is an odd integer.

If $F$ is surjective onto $\mathbb{N}$, then $X_n$ has infinitely many integral points (for each $z \in \mathbb{N}$, $z^n$ is in the image of $F$, and so $z^n = F(x, y)$ for some $x, y \in \mathbb{Z}$). We may therefore assume $X_n$ has infinitely many integral points. If these points are not dense on $X_n$, they must be concentrated in finitely many curves. As punctured curves of genus $g \geq 1$ only have finitely many integral points, all but finitely many points lie on rational curves. By replacing $n$ by a multiple of itself, these rational curves may be replaced by cyclic covers. 

\begin{lemma}
Let $f$ be a function on $\mathbb{P}^1$, then the curve $s^n = f(t)$ is irreducible and has positive genus for some suitable $n$, unless $f$ only has one zero and one pole.
\end{lemma}
\begin{proof}
Suppose $f$ has two or more zeroes or poles. Let $n \geq 3$ be a prime which does not divide the order of the zeroes and poles, then the curve is irreducible and $s^n = f(t)$ ramifies over the zeroes and poles. By Riemann-Hurwitz, the genus of the cover, $g$, satisfies $2g - 2 = n(-2) + \sum (e_P - 1) \geq -2n + 3(n - 1) = n - 3 \geq 0$.
\end{proof}

Therefore, after taking $n$ large enough, it may be assumed that when $F$ is restricted to any of these rational curves with dense points it has only one zero and one pole. Fix one of these curves and its image, $C$, after applying the projection to $(x, y)$. By suitable parameterisation of $C$, $\phi : \mathbb{P}^1 \dashrightarrow C$, it can be assumed that $(F \circ \phi)(t) = \frac{1}{A}t^{\ell}$ for some integer $A$, and a positive integer $\ell$. As integral values of $F$ are restricted to integer $t$, the functions $x(\phi(t))$ and $y(\phi(t))$ must be integral infinitely often for integer input. The following lemma shows $x(\phi(t))$ and $y(\phi(t))$ are polynomials as a result.

\begin{lemma}
Let $\frac{f}{g}$ be a rational function, with $f$ and $g$ coprime. Suppose that for infinitely many integers, $t$, $\frac{f(t)}{g(t)} \in \mathbb{Z}$, then $g$ is constant. 
\end{lemma}
\begin{proof}
By scaling, $f$ and $g$ may be assumed to be integral. By coprimality, there exists $r, s \in \mathbb{Z}[T]$ and $N \in \mathbb{Z}$ such that $rf + sg = N$. Dividing by $g$ shows $r \frac{f}{g} + s = \frac{N}{g}$, and the left hand side is integral infinitely often. As $g$ is always integral, this shows $g(x) | N$ for infinitely many $x$, but $N$ has only finitely many divisors, and so $g$ must take the same value infinitely often, and so is constant.
\end{proof}

If the image of $F$ is contained in $\mathbb{N}$, then $\ell$ has to be even, otherwise $\frac{1}{A} t^{\ell}$ would be negative for $t < 0$ (or $t > 0$ if $A < 0$), but infinitely many integer $(x, y)$ come from $t < 0$ (or $t > 0$), since $x$ and $y$ are polynomials which are integer infinitely often for $t > 0$ (or $t < 0$). It follows that the perfect $n$th powers represented by $F$ are all squares, which are a zero density set in the $n$th powers. As there are only finitely many curves to consider, $F$ cannot represent all $n$th powers, for $n$ odd, and so is not surjective. This proves Theorem~\ref{fewpts}.

\section{Geometry of $X_n$ and Complements of Plane Curves}
In this section we study the geometry of cyclic covers of the affine plane.

As before, $X_n$ is the surface $z^n = F(x, y)$. By shifting $F$ by a constant, we can assume the curve $C : F(x, y) = 0$ is smooth in $\mathbb{A}^2$, this either adds finitely many negative values into the image of $F$ or removes finitely many positive ones. As $C$ is smooth, this makes $X_n$ smooth. 

\begin{lemma}
The surface $X_n$ is of log general type for large enough $n$, if $\mathbb{A}^2 \setminus \{F(x, y) = 0\}$ is.
\end{lemma}
\begin{proof}
As for canonical divisors of projective varieties, there is a Riemann-Hurwitz formula for logarithmic canonical divisors of open surfaces \cite{Iitaka}. As $\pi_n : X_n \to \mathbb{A}^2$ is finite and ramified only over $F(x, y) = 0$, we have $\overline{K}_{X_n} \sim \pi_n^* \overline{K}_{\mathbb{A}^2} + (n - 1) D$, where $D$ is the locus of $z = 0$ in $X_n$. As $nD = \pi_n^*  D'$, where $D'$ is the locus $F(x, y) = 0$, $n \overline{K}_{X_n} \sim \pi_n^* (n \overline{K}_{\mathbb{A}^2} + (n - 1) D')$. The Iitaka dimension of $n \overline{K}_{\mathbb{A}^2} + (n - 1) D'$, for large $n$, is, definitionally, $\underline{\kappa}_{\mathbb{A}^2 \setminus D'}$. If the logarithmic Kodaira dimension, $\overline{\kappa}$, is 2 for any surface, then $\underline{\kappa} = 2$ by a result of Iitaka \cite{Iitaka}. Therefore, if $\mathbb{A}^2 \setminus D'$ is of log general type, then $n \overline{K}_{X_n}$ is big for large enough $n$, and the result follows.
\end{proof}

Combining this lemma with Theorem 1 and Vojta's conjecture gives Proposition 1.

The geometry of curves $C$ such that $\mathbb{A}^2 \setminus C$ is not of log general type is determined in Aoki's ``\'{E}tale Endomorphisms of Smooth Affine Surfaces'' \cite{aoki}. They prove the following classification of such curves, as the vanishing of $H$, over $\overline{\mathbb{Q}}$.
\begin{theorem}
Let $C : H(x, y) = 0$ be a curve such that $\kappa(\mathbb{A}^2 \setminus C) \neq 2$, then $H(x, y) = a$ is rational for all but finitely many $a$, and moreover, one of the following holds:
\begin{enumerate}
 \item The curve $C$ is isomorphic to $\mathbb{A}^1$.
 \item The curve $C$ is isomorphic to $\mathbb{A}^1_* = \mathbb{A}^1 \setminus \{0\}$.
 \item There exists a polynomial $G \in \overline{\mathbb{Q}}[x, y]$ such that $V(G)$ is isomorphic to $\mathbb{A}^1$, and $C$ meets $V(G)$ at a single point, and the same is true for $V(G - a)$ for general $a$.
 \item The curve $C$ is singular, and after a suitable change of variables, it is of the form $y^n = x^m$ for coprime $n, m$.
\end{enumerate}
\label{logkod}
\end{theorem}

We apply this theorem to $F(x, y) = 0$. By assumption, the vanishing of $F$ is assumed to be smooth, so case 4 can be automatically ruled out.

To understand the remaining cases, the following theorem of Abhyankar and Moh \cite{abymoh} is helpful.
\begin{theorem}
Let $\phi : \mathbb{A}^1 \to \mathbb{A}^2$ be an embedding, and write $\phi(t) = (f(t), g(t))$, where $f, g \in K[t]$, for any characteristic 0 field $K$ (not necessarily algebraically closed), then either $\mathrm{deg}(f) \mid \mathrm{deg}(g)$ or $\mathrm{deg}(g) \mid \mathrm{deg}(f)$
\label{abmoh}
\end{theorem}

There is a geometric corollary of this statement that is often more convenient, and we use similar methods to the proof given in \cite{abymoh}.

\begin{corollary}
Let $X \subset \mathbb{A}^2$ be geometrically smooth, defined over $K$ (any characteristic 0 field), and isomorphic to $\mathbb{A}^1$ then there is a change of variables, also defined over $K$, such that $X$ is a coordinate axis.
\label{abmohcor}
\end{corollary}
\begin{proof}
As $X$ is isomorphic to $\mathbb{A}^1$, it has a parameterisation $\phi = (f, g) : \mathbb{A}^1 \to \mathbb{A}^2$. By Theorem~\ref{abmoh}, either $\mathrm{deg}(f) \mid \mathrm{deg}(g)$ or $\mathrm{deg}(g) \mid \mathrm{deg}(f)$. If, say, $\mathrm{deg}(g) = n \times \mathrm{deg}(f)$, then by applying the change of variables $(x, y) \mapsto (x, y - ax^n)$, for some suitable choice of $a$, the degree of $g$ decreases. This can be repeated until the degree of $g$ is smaller than $f$, and then the roles reverse. Eventually, one degree is 0, and the other is necessarily 1 (as the curve is reduced). After translation, the curve is then an axis.
\end{proof}

We also need the classification of $f$ satisfying condition 2, which gives an anologous statement to the geometric version of Abhyankar and Moh in this case. This is given in \cite{aoki}, from earlier work of Miyanishi and Sugie \cite{miysug}.
\begin{theorem}
Suppose $H \in \overline{\mathbb{Q}}[x, y]$ is generically rational (i.e. $H(x, y) = a$ is rational for all but finitely many $a$), and the general curve $H = a$ is missing 2 points from its projective closure, then there is a change of variables such that $f$ is, up to scaling, of the form $x^a y^b + 1$ or $x^a (x^\ell y + p(x))^b + 1$, where $a, b, \ell > 0$, $(a, b) = 1$, and $p$ is a polynomial of degree less than $\ell$ with $p(0) \neq 0$.
\end{theorem}

The following lemma controls the asymptotic behaviour of such rational curves, in a way analogous to Theorem~\ref{abmoh}.

\begin{lemma}
Let $\phi : \mathbb{A}_*^1 \to \mathbb{A}^2$ be a parameterisation of a rational curve, $H(x, y) = 0$, where $H$ is generically rational, with $\phi(t) = (f(t), g(t))$ for some $f, g \in K[t, t^{-1}]$. Then either one of $f$ or $g$ does not have a pole at infinity, or the order of the pole of $f$ at infinity divides that of $g$, or vice versa. 
\label{growth}
\end{lemma}
\begin{proof}
Suppose both $f$ and $g$ have a pole at infinity, and the order of the pole of $f$ does not divide that of $g$ and the order of the pole of $g$ does not divide the order of the pole of $f$. There is a change of variables which puts it into the form as in the preceding theorem. In those coordinates, the parameterisation is $(r(t), s(t))$, where one of $r$ and $s$ does not have a pole at infinity, depending on which of the two missing points $\infty$ is mapped to. The change of variables must therefore decrease the order of the pole at infinity for at least one coordinate.

Define a partial ordering, $\preccurlyeq$, on $\mathbb{N}^{(2)}$ by $(a, b) \prec (c, d)$ if $\min(a, b) \leq \min(c, d)$ and $\max(a, b) \leq \max(c, d)$. All changes of variable can be obtained by performing the following elementary operations: a linear transformation (only scaling $x$ and $y$ or swapping them); adding a polynomial in $x$ to $y$; adding a polynomial in $y$ to $x$. Let $(a_1, b_1), ..., (a_n, b_n)$ be the degrees of $f$ and $g$ during a minimal sequence of elementary operations such that $(a_n, b_n) \preccurlyeq (a_{n - 1}, b_{n - 1})$, and, moreover, there is no other sequence of transformations with the same end points (in terms of coordinates) such that $(a_1, b_1) = (a_1', b_1') \preccurlyeq ... \preccurlyeq (a_m', b_m') = (a_n, b_n)$. 

As $(a_n, b_n) \preccurlyeq (a_{n - 1}, b_{n - 1})$, the last transformation has to be adding a power of one variable to another. We assume it is adding a polynomial in $x$ to $y$, the other case is proved similarly. As this operation will not give a reduction in degree when starting with $f$ and $g$ (as neither degree divides the other), the sequence must contain at least 2 tranformations.

If the penultimate operation is a linear transformation, it is either a scaling or it interchanges $x$ and $y$. In the former case, this scaling can be incorporated into the polynomial giving a shorter sequence. If it swaps $x$ and $y$, then swapping the role of $x$ and $y$ in the last operation and removing the swap gives a shorter sequence. 

If instead, the penultimate step is adding a polynomial in $x$ to $y$, then this can be combined with the last operation to give a shorter sequence or doing so means the degrees do not decrease (and so the sequence can be replaced by one with constantly increasing degree).

This leaves only adding a polynomial in $y$ to $x$. This makes $a_{n - 1} \geq b_{n - 1}$, unless there is cancellation, which would make $a_{n - 1} < a_{n - 2}$, $b_{n - 1} = b_{n - 2}$, giving a shorter sequence. The last step must then be subtracting a multiple of $x$ from $y$, and the penultimate step adding a multiple of $y$ to $x$. The combined map is $x \mapsto x + ay$, $y \mapsto (1 + ab)y + bx$ for some scalars $a, b$. If $1 + ab \neq 0$, the degree reduction in $y$ can be achieved by just adding a suitable multiple of $x$, giving a shorter sequence. Otherwise, the map is $x \mapsto x + ay$, $y \mapsto bx$. For this to give an inequality in degrees in the final step, $a_{n - 2} = b_n < b_{n - 1} = b_{n - 2}$, but then the sequence of transformations can be written without the inequality at the final step by swapping $x$ and $y$, and scaling, (degree $(b_{n - 2}, a_{n - 2})$) before adding a multiple of what is now $y$ to $x$, which as $a_{n - 2} < b_{n - 2}$ does not change the degrees.
\end{proof}

\section{Arithmetic of Generically Rational Polynomials}
In this section, we will show that the polynomials arising in the classification of Theorem~\ref{logkod} either fail to represent all positive integers, or represent arbitrarily large negative integers. An important result is the following theorem, communicated to the author by Samir Siksek (personal communication, 2023).

\begin{theorem}
Let $C$ be a geometrically rational smooth curve over $\mathbb{Q}(T)$ such that for all but finitely many $t \in \mathbb{N}$, the specialisation of $C$ at $T = t$ has a rational smooth point, then $C$ is isomorphic to $\mathbb{P}^1 / \mathbb{Q}(T)$.
\label{hasse}
\end{theorem}
\begin{proof}
By taking the anti-canonical embedding, $C$ can be assumed to be a plane conic, and moreover its equation is $f(T)X^2 + g(T)Y^2 = Z^2$ for some non-zero rational functions $f, g$. If $C$ has a smooth point after specialising at $t$, then so does the conic, and the Hilbert symbol $(f(t), g(t))$ vanishes. As this holds for all but finitely many $t \in \mathbb{N}$, the Hilbert symbol $(f(T), g(T))$ also vanishes by a theorem of Serre \cite{serre}. This shows that the conic is isomorphic to $\mathbb{P}^1 / \mathbb{Q}(T)$, and so $C$ is as well.
\end{proof}

Applying the theorem to the curve $F(x, y) = T$ shows that it is a rational curve over $\mathbb{Q}(T)$. In particular, by taking the completion at the infinite place of $\mathbb{Q}$, $F : \mathbb{R} \times \mathbb{R} \to \mathbb{R}$ is surjective. 

The other main tool is this proposition which guarantees the existence of certain integer points, even after rational change of variable, that is, after applying a sequence of the following: linear transformations (from $\mathrm{GL}_2(\mathbb{Q})$); adding a rational polynomial in $x$ to $y$; and adding a rational polynomial in $y$ to $x$.

\begin{propn}
Suppose $R \subset \mathbb{Q}^2$ satisfies weak approximation in the following sense. For every $(x_0, y_0) \in \mathbb{Q}^2$, $N \in \mathbb{N}$, there exists an $(x, y) \in R$ such that $x - x_0, y - y_0 \in N \mathbb{Z}$. Such a set contains infinitely many integer points, and moreover, its image under any rational change of variables also satisfies the hypothesis.
\label{wkapp}
\end{propn}
\begin{proof}
Let $(x_0, y_0)$ be an integer point not in $R$ (if all integer points are in $R$, then $R$ contains infinitely many integer points), then by applying the hypothesis for increasing $N$ gives infinitely many integer points in $R$.

As any change of variables can be expressed as a combination of linear transformations, adding powers of $x$ to $y$, and adding powers of $y$ to $x$, it is enough to check invariance under each type. In combination with the other 2 types of transformation, linear transformations can be assumed to be either scalings, or swapping the two variables.

The weak approximation property is preserved by swapping $x$ and $y$. We can also assume the scaling affects only one variable, and by symmetry that is $x$. Let the scale factor be $\alpha = \frac{a}{b}$ where $a \in \mathbb{Z}$, $b \in \mathbb{N}$ are coprime. Given $(x_0, y_0) \in \mathbb{Q}^2$ and $N \in \mathbb{N}$, there exists $(x, y) \in R$ such that $x - \alpha^{-1} x_0, y - y_0 \in bN\mathbb{Z}$, therefore $(\alpha x, y)$ is in the image of $R$ and $ax - bx_0 \in abN \mathbb{Z}$, so $\alpha x - x_0 \in N \mathbb{Z}$. 

By symmetry, it is again enough to check one of the other two cases. We check the case of $(x, y) \mapsto (x, y + x^n)$. If $n = 0$, then it clearly preserves this property, so take $n > 0$. For any pair $(x_0, y_0) \in \mathbb{Q}^2$ and $N \in \mathbb{N}$, let $x_0 = \frac{a}{b}$, then there exists $(x, y) \in R$ such that $x - x_0, y - y_0 + x_0^n \in b^{n - 1} N \mathbb{Z}$. The point $(x, y + x^n)$ is in the image of $R$ and $(y + x^n) - y_0 = (y - y_0 + x_0^n) + (x^n - x_0^n)$. The first term is in $N\mathbb{Z}$ by assumption, the second term factors as $(x - x_0)(x^{n - 1} \ldots + x_0^{n - 1})$, the first term is integer and divisible by $b^{n - 1} N$, the second has denominators at worst $b^{n - 1}$, and so their product is an integer divisible by $N$. 
\end{proof}

Our strategy to show no generically rational polynomial, $F$, gives a surjection $\mathbb{Z}^2 \to \mathbb{N}$ is the following: for each case in the classification of Theorem~\ref{logkod} use rational changes of variable to simplify the form of the polynomial to $G$, and use that $F$, and so $G$, surjects onto $\mathbb{R}$ to conclude that $\{(x, y) \mid G(x, y) < 0\}$ satisfies weak approximation, as in the proposition, and so the same is true for $F$. Not only then does $F$ have infintely many integer points where it is negative, but these integer points represent arbitarily large negative numbers since shifts of $F$ still meet the conditions of the theorem. We start under the assumption of condition 1.

\begin{theorem}
Functions, $F$, such that $F(x, y) = a$ is generically isomorphic to $\mathbb{A}^1$, represent arbitrarily large negative integers, as $\{(x, y) \mid F(x, y) < 0\}$ satisfies weak approximation.
\label{case1}
\end{theorem}
\begin{proof}
By the theorem of Abhyankar-Moh, Theorem~\ref{abmoh}, there is a rational change of variables to $(u, v)$ such that $F(u, v) = u$. Here, the set where $F$ is negative is $u < 0$. This satisfies weak approximation, if $(u_0, v_0) \in \mathbb{Q}^2$ has $u_0 < 0$, then for any $N$, $(u_0, v_0)$ fulfills the hypotheses, and if $u_0 \geq 0$, then some multiple of $N$, say $Nd$, is larger than $u_0$, and so $(u_0 - Nd, v_0)$ is in the set. 
\end{proof}

The next case to consider is the second, where $C \simeq \mathbb{A}^1_*$. We start by observing the two removed points must be defined over $\mathbb{R}$. Suppose not, then the parameterisation, $(f, g) = \phi : \mathbb{P}^1\setminus \{a, b\} \to \mathbb{A}^2$, of $F(x, y) = 0$ has no poles in $\mathbb{R}$. In particular, $\phi(\mathbb{R})$ is bounded and connected. The region bounded by this curve is therefore either the region where $F$ takes positive values or negative values. As this region is closed and bounded, it is compact, and $F$ has bounded image on it. Either it fails to represent arbitrarily large positive integers and so $F(\mathbb{Z} \times \mathbb{Z}) \neq \mathbb{N}$, or it fails to represent large negative reals, contradicting the remarks following Theorem~\ref{hasse}. We therefore assume from now on that both $a$ and $b$ are real.

After compactifying $\mathbb{A}^2$ to $\mathbb{P}^2$, the parameterisation $\phi$ extends to $\tilde{\phi} : \mathbb{P}^1 \to \mathbb{P}^2$.

\begin{theorem}
Let $\phi$ and $C$ be as above. Assume further that the images of the two deleted points are distinct in $\mathbb{P}^2$, then the set $\{(x, y) \mid F(x, y) < 0\}$ satisfies weak approximation.
\label{case2a}
\end{theorem}
\begin{proof}
We proceed by 2 cases.
\paragraph{Case 1:} We assume there exists a line, $L \subset \mathbb{P}^2$, such that $(0 : 0 : 1) \in L$, $\tilde{\phi}(a), \tilde{\phi}(b) \notin L$, and $F(L \cap \mathbb{A}^2(\mathbb{R}))$ is not bounded below. 
Let $L_1$, $L_2$ be the lines such that $(0 : 0 : 1) \in L_1, L_2$ and $\tilde{\phi}(a) \in L_1$, $\tilde{\phi}(b) \in L_2$. Further, choose lines $L_3$ and $L_4$, such that $L_3$ is between $L_1$ and $L$, and $L_4$ is between $L_2$ and $L$. We now restrict to $\mathbb{A}^2(\mathbb{R}) \subset \mathbb{P}^2$. As $F$ is not bounded below on the real points of $L$, there is a half-line of $L$ for which $F < 0$. As $F$ only changes sign along $C$, this shows there is a non-compact region, $S$, bounded by $C$, containing the half-line of $L$, such that $F < 0$ on the interior. As $L_3$ is between $L$ and $L_1$, and $L_4$ is between $L$ and $L_2$, $C$ is bounded on the region, $T$, bounded by $L_3$ and $L_4$. In particular, inside $S \cap T$, there is a line such that $C$ and $(0, 0)$ lie on the same side of the line. Let $R$ be the intersection of $T$, and the other side of this line than $(0, 0)$. An example of this process is shown in Figure~\ref{constr}, with $R$ shaded.

As $R$ is contained in both $T$, and is on the correct side of the line bounding $C$, $R$ is contained in $S$, so $F$ is negative in $R$. The 3 lines defining $R$ are not parallel, and do not define a bounded region, so $R$ satisfies weak approximation.

\begin{figure}
\begin{tikzpicture};
\fill[lightgray] (1.77, 0.89) -- (1.77, 3.58) -- (2.5, 5) -- (5, 5) -- (5, 2.5);
\draw (-2, 0) -- (5, 0) node[right] {$L_1$};
\draw (0, -1) -- (0, 5) node[above] {$L_2$};
\draw[domain=0.56:5, samples=40] plot (\x, {(\x + 1)/(\x^2)}) node[above] {$C$};
\draw[domain=0.36:2, samples=20] plot (-\x, {(-\x + 1)/(\x^2)});
\draw (-1, -1) -- (5, 5) node[right] {$L$};
\draw[dashed] (-0.5, -1) -- (2.5, 5) node[right] {$L_4$};
\draw[dashed] (-2, -1) -- (5, 2.5) node[right] {$L_3$};
\end{tikzpicture}
\caption{An illustration of the region used in case 1 of the proof of Theorem~\ref{case2a}.}
\label{constr}
\end{figure}
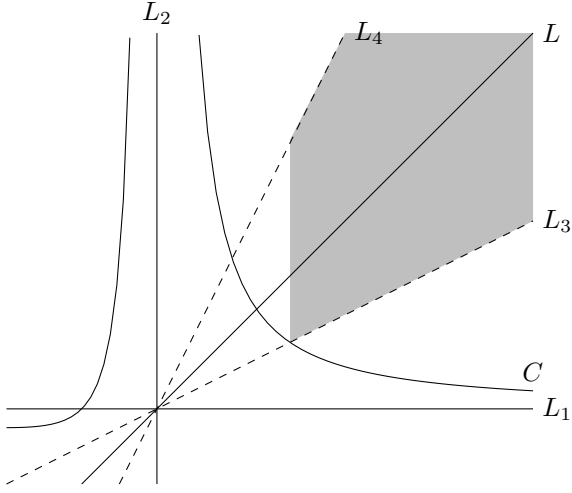

\paragraph{Case 2:} Suppose instead, no such line exists. Take 2 lines, $L_1$, $L_2$, containing $(0 : 0 : 1)$, so that $\mathbb{P}^2 \setminus (L_1 \cup L_2)$ is the union of two disjoint open sets, $U_1$, $U_2$ with $\tilde{\phi}(a) \in U_1$ and $\tilde{\phi}(b) \in U_2$. By assumption, $F$ is bounded below on these two lines, say by $X$. Let $c < X$ be sufficiently general, so that $F(x, y) = c$ is smooth, rational, and is isomorphic to $\mathbb{A}^1_*$. As $c < X$, $F(x, y) = c$ is disjoint from $L_1$ and $L_2$, and so in particular, lies within $U_1 \cup U_2$. As $F(x, y) = c$ is rational, its closure in $\mathbb{P}^2$ is connected. However, its closure must pass through both $\tilde{\phi}(a)$ and $\tilde{\phi}(b)$, which would make $U_1 \cup U_2$ path-connected, which is a contradiction.

\end{proof}

We can now understand curves with a single point at infinity, which corresponds to 2 points in the normalisation, by reducing to the case where the points are separate, except in one exceptional circumstance.

\begin{theorem}
Let $\phi$ and $C$ be as before, and assume the images of the deleted points coincide in $\mathbb{P}^2$, then $\{F(x, y) < 0\}$ satisfies weak approximation.
\label{case2b}
\end{theorem}
\begin{proof}
As there is a unique point at infinity, it must be rational, and so by a rational change of coordinates it can be assumed to be $(0 : 1 : 0)$. By Lemma~\ref{growth}, either $f$ or $g$ does not have a pole at $a$ or the order of the two poles divide each other. In the former, it must be $f$ that does not have a pole, otherwise it would not go through $(0 : 1 : 0)$. The same statement holds at $b$ as well. If in both cases one function has no pole, then as it must be $f$ in each case, $f$ has no poles. This would make $f$ a constant, and as $g$ has two poles, this would be a double cover of a line, and so can be discarded. 

We can now assume that both $f$ and $g$ have a pole at, at least one of, $a$ and $b$, say $a$. As the degree of $f$ divides the degree of $g$ (the degree of $g$ must be larger as the curve goes through $(0 : 1 : 0)$ at $a$), the asymptotic behaviour of $C$ is $y \sim A x^d$, where $d$ is integer. If $A$ is rational, then the change of variables $(x, y) \mapsto (x, y + A x^d)$ is defined over $\mathbb{Q}$. This change of variable decreases the order of the pole of $g$ at $a$. This process can be repeated until one of the following holds (as the degrees of $f$ and $g$ decrease at every step): the points at infinity are separate, the curve no longer passes through $(0 : 1 : 0)$, or $A$ is no longer rational. Let $G$ be the result of composing $F$ with the changes of variable, and let $u$ and $v$ be the new coordinates.

In the first case, by the preceding theorem combined with Proposition~\ref{wkapp}, $\{(x, y) \mid F(x, y) < 0\}$ satisfies weak approximation. In the second, either the point at infinity is now $(1 : 0 : 0)$, where by swapping $u$ and $v$ it can be returned to $(0 : 1 : 0)$ without increasing the degrees of $f$ and $g$, or instead, subtracting a further multiple of $u$ from $v$ will move the point at infinity to $(1 : 0 : 0)$, further decreasing the order of $f$ and $g$.

The final case is when $A$ is not rational. As $f$ and $g$ are defined over $\mathbb{Q}$, this shows that $a$ and $b$ are not, and so are defined over a (real) quadratic extension, $\mathbb{Q}(\sqrt{D})$ In this case, let $A = \alpha + \beta \sqrt{D}$, then the growth at $b$ must be given by $v \sim (\alpha - \beta \sqrt{D}) u^d$. Subtracting $\alpha u^d$ from $v$ makes the two asymptotics have opposite signs. This means for large enough $u$, either $G(u, v) < 0$ if $v > C_1 u^d$ or $v < -C_2 u^d$ for suitable $C_1, C_2 > 0$; or if $-C_3 u^d < v < C_4 u^d$for some $C_3, C_4 > 0$. Both of these sets satisfy weak approximation, therefore by Proposition~\ref{wkapp}, $\{(x, y) \mid F(x, y) < 0\}$ satisfies weak approximation.
\end{proof}

Combining the two theorems shows that for any $F$ satisfiying condition 2 of Theorem~\ref{logkod}, either $F$ does not represent (over $\mathbb{R}$) all real numbers, or $\{F(x, y) < 0\}$ satisfies weak approximation.

The final case is those satisfying condition 3 of Theorem~\ref{logkod}: there exists a polynomial $G \in \overline{\mathbb{Q}}[x, y]$ such that $V(G - a)$ is rational for general $a$, and these curves meet $C$ in one point (depending on $a$). We start by showing $G$ can be assumed to be defined over $\mathbb{Q}$.

\begin{propn}
The auxilliary polynomial $G$, after a linear transformation, is defined over $\mathbb{Q}$, or $F$ also satisfies conditions 1 or 2.
\end{propn}
\begin{proof}
We assume that $F$ does not satisfy conditions 1 or 2, that is, $C$ is missing at least 3 points.

Over the field of definition of $G$, there is a change of variables so that $G$ becomes a coordinate by applying Corollary~\ref{abmohcor}. Let the new coordinates be $u, v$, with $G(u, v) = u$. By the intersection assumption, $F$ must be given by $f_1(u)v + f_2(u)$ for some polynomials $f_1, f_2$, and by shifting $v$ by a power of $u$, the degree of $f_2$ can be assumed to be less than the degree of $f_1$. If $F$ does not fit into any of the other classifications, $f_1$ has at least 2 roots (otherwise $C$ is isomorphic to $\mathbb{P}^1$ missing one or two points). Suppose $G(x, y)$ is not defined over $\mathbb{Q}$, then there exists a conjugate of $G$, $G'$. The vanishing of $G'$ must be a rational curve meeting the transformed version of $C$ in one place. Let the parameterisation of $V(G')$ be $(r(t), s(t))$, then $f_1(r(t)) s(t) + f_2(r(t))$ has a single root. Since replacing the vanishing of $G'$ with the vanishing of a shift of $G'$ does not change the number of intersections, $f_1(r(t)) s(t) + f_2(r(t))$ is linear. If the degree of $r$ is at least 1, then the degree of $f_1(r(t))$ is larger than $f_2(r(t))$, and the degree of $f_1(r(t))$ is at least 2 (as $f_1$ has at least 2 roots). This precludes $f_1(r(t)) s(t) + f_2(r(t))$ being degree 1, and so $r$ is a constant. This makes $G'(u, v) = au + b = aG + b$ for some $a, b$. This equality must, therefore, hold for $G(x, y)$ and $G'(x, y)$. Suppose the coefficient of $x^i y^j$ in $G$ is non-zero, where at least one of $i$ and $j$ is non-zero. As scaling $G$ does not change any of the properties, this coefficient can be assumed to be 1. By comparing the coefficients of $x^i y^j$ in $G(x, y)$ and $G'(x, y)$, it follows that $a = 1$, so $G' = G + b$. This shows $G$ can be written as a rational polynomial, plus an algebraic number, $c$. Shifting $G$ by a constant does not change the properties of $G$, so $c$ can be assumed to be 0, and $G$ is now defined over $\mathbb{Q}$.
\end{proof}

With this proposition, this case proceeds almost as in Theorem~\ref{case1}.
\begin{theorem}
Suppose $F$ meets condition 3, and not any other, then $\{(x, y) \mid F(x, y) < 0\}$ satisfies weak approximation.
\label{case3}
\end{theorem}
\begin{proof}
By the previous proposition, we can replace $G$ by one with the same properties defined over $\mathbb{Q}$, and so there exists a rational change of variables so that $F(x, y)$ becomes $H(u, v) = f_1(u)v + f_2(u)$ with the degree of $f_2$ smaller than the degree of $f_1$. The set $\{(u, v) \mid H(u, v) < 0\}$ satisfies weak approximation, since for any $(u_0, v_0) \in \mathbb{Q}^2$ and $N \in \mathbb{N}$, there is a shift $u_0 + Nd_1$ such that $f_1(u_0 + Nd) \neq 0$, and a shift $v_0 + Nd_2$ such that $f_1(u_0 + Nd)(v_0 + Nd_2) + f_2(u_0 + Nd) < 0$. By Proposition~\ref{wkapp}, this shows the same is true for $\{(x, y) \mid F(x, y) < 0\}$.
\end{proof}

Combining Theorems~\ref{case1}, \ref{case2a}, \ref{case2b} and \ref{case3} exhausts the conditions of Theorem~\ref{logkod}. Any polynomial $F(x, y)$ such that $z^n = F(x, y)$ has dense integral points for all odd $n$ must therefore either not represent all positive integers (as it does not represent large negative reals, using Theorem~\ref{hasse}) or it represents large negative integers.
\printbibliography
\end{document}